\documentclass{article}
\usepackage{graphicx} 
\usepackage{algpseudocode}
\usepackage{amsmath}
\usepackage{algorithm}
\title{Adaptively Learning Memory Incorporating PSO}

\date{February 2024}

\begin{document}
\maketitle

\author{Dmytro Shchyrba}, drumyraetherid@gmail.com \author[1]{Izabela Aleksandra Paniczek}\author[2]
Wroclaw University of Science and Technology, Wroclaw, Poland

\section{Abstract}

\textwidth 14.5cm

\setcounter{page}{1}

\renewcommand{\labelenumii}{\arabic{enumi}.\arabic{enumii}}
\providecommand{\keywords}[1]
{
  \small	
  \noindent
  \textbf{Keywords:} #1
}

\noindent Selection of perefect parameters for low-pass filters can sometimes be an expensive problem with no analytical solution or differentiability of cost function. In this paper, we introduce a new PSO-inspired algorithm, that incorporates the positive experiences of the swarm to learn the geometry of the search space,thus obtaining the ability to consistently reach global optimum and is especially suitable for nonsmooth semiconvex functions optimization. We compare it to a set of other algorithms on test functions of choice to prove it's suitability to a certain range of problems, and then apply it to the problem of finding perfect parameters for exponential smoothing algorithm.

\keywords{Swarm intelligence, Time Series Analysis, Intelligent Systems}

\section{Introduction}

\noindent Optimization problems are fairly common in most complex human endeavors, such as engineering \cite{rajput2020review} \cite{nutakki2023review}, logistics \cite{zhou2023application} or artificial intelligence \cite{abdulkadirov2023survey}. Currently, there are lots of methods available for finding either optimal or sub-optimal solutions for them, such as nature-based and gradient-based methods. Although gradient-based methods tend to have faster convergence to a solution, requirement of Jacobian and sometimes even Hessian matrix, becomes a strong disadvantage in real-life applications due to highly complicated objective functions, usually making computation of the derivatives unfeasible or even physically impossible. Thus, although perceived as slower comparatively to gradient-based, population-based algorithms are still sufficient for a wide variety of tasks. 
Our work is mostly interested in nonsmooth semiconvex nondifferentiable functions, such as loss function in choosing the perfect parameters for filters such as exponential smoothing - where the number of cost evaluations is supposed to be as low as possible, but the cost function itself is nondeterministic and there is neither an analytical solution, nor possibility to apply gradient-based algorithms.

\section{Backgrounds on PSO}
\noindent Particle Swarm Optimization \cite{kennedy1995particle} is one of the most popular optimization algorithms.Being inspired by the social behaviour of birds, it has proven itself worthy in many cases and application. In PSO, search agents possess 2 cognitive factors for searching- social and individual.Each iteration,the coordinates are updated according to the equations
    \begin{equation}
        x_i (t+1)=x_i (t)+v_i (t+1) 
    \end{equation}
    \begin{equation}
        v_i (t+1)= w\cdot r_i \cdot v_i(t) + c_1\cdot r_i\cdot(g_i (t) - x_i (t)) + c_2\cdot r_i\cdot(p_i (t) - x_i (t))
    \end{equation}
Where $x_i$ is the coordinates vector for i-th particle, $v_i$ is respectively the velocity, $w$ is the inertia hyperparameter, $r_i$ is the random number generated from uniform distribution between, and $c_1$ and $c_2$ are the parameters for  social and indiviual cognitive factors.

\section{Preliminaries on Exponential Smoothing}

Exponential smoothing is a widely used technique for making short-term forecasts in time series analysis. Its popularity stems from the simplicity, efficiency, and flexibility it offers, making it suitable for a broad range of applications. The core principle behind exponential smoothing is to assign exponentially decreasing weights to past observations, with the most recent observations given more significance.

\subsection{Basic Concept}

The simplest form of exponential smoothing, often referred to as Simple Exponential Smoothing (SES), can be expressed as:
\begin{equation}
    S_t = \alpha x_t + (1 - \alpha) S_{t-1}
\end{equation}
where:
\begin{itemize}
    \item $S_t$ is the smoothed statistic for the current period,
    \item $x_t$ is the actual observation at time $t$,
    \item $S_{t-1}$ is the smoothed statistic for the previous period,
    \item $\alpha$ is the smoothing constant, a value between 0 and 1.
\end{itemize}

\subsection{Importance of Smoothing Constant $\alpha$}

The smoothing constant, $\alpha$, plays a critical role in the effectiveness of exponential smoothing models. It determines the rate at which the influence of older observations declines. A higher value of $\alpha$ places more emphasis on recent observations, making the forecast more responsive to changes. Conversely, a lower value of $\alpha$ results in a smoother forecast, which might be desirable in the presence of significant random fluctuations in the data.

\subsection{Extensions and Variations}

Exponential smoothing has been extended to handle more complex time series patterns, such as trends and seasonality. The two most notable extensions are:
\begin{itemize}
    \item Holt’s Linear Trend Method: Introduces an additional parameter to capture the trend in the data.
    \item Holt-Winters Seasonal Method: Incorporates seasonal components into the forecast, handling both trends and seasonality.
\end{itemize}

\subsection{Parameter Optimization}

Finding the optimal parameters for an exponential smoothing model, including the smoothing constants for the level, trend, and seasonal components, is crucial for achieving accurate forecasts. Various optimization techniques can be employed to minimize the forecast error, including grid search, gradient descent, and evolutionary algorithms. The choice of optimization method and error metric (e.g., Mean Squared Error, Mean Absolute Error) can significantly impact the model's performance.

\subsection{Conclusion}

Exponential smoothing offers a robust framework for time series forecasting. By carefully selecting and optimizing its parameters, practitioners can adapt the method to a wide range of forecasting scenarios, balancing responsiveness against the desire for smoothness in their forecasts.

\section{Proposed approach}
The algorithm we are proposing consists of 3 main parts: incorporation of the preserved history into the swarm thinking abilities, reprojection across the best neighbors and natural selection component

\begin{enumerate}
\item \underline{Applications}
The algorithm is especially focused on the multimodal large scale optimization tasks with a pattern towards the global optima.As it is shown in results section,it drastically outperforms it's opponents on the functions where optimization landscape has a tendency to bend towards the global optima.To obtain such abilities,three factors have been considered. 
Firstly, the direction choice. Inspired by the Unified PSO, after usual velocity update, particles are being projected across their best-performing neighbors.
Secondly,the consistency of the learned tendencies of landscape being optimized
In order to make the algorithm less reliant on manually chosen parameters and free it from the demanded prior knowledge, it is learned on the run from preserved history, concept of which will be explained in more detail in the next section.
\item \underline{Preserved history} \\ To incorporate the information from the previous generations in a more effective way, we store the information about the best fitness values and their respective solutions in a set of size $N_p$, which is a custom parameter. Each of the objects which it contains possess 3 attributes-the coordinates, the fitness value, and the ,,age’’ (the iteration when the solution was appended to the set).
\item \underline{Reprojection} \\ The proposed approach involves two kinds or reprojections – involving individual and collective thinking components
\begin{enumerate}
    \item \underline{Individual} \\ Here is the update equation for an individual particle in
    \begin{equation}
        x_i (t+1)=x_i (t)+w_i (t+1)+2\cdot(l_i (t)-x_i (t)) 
    \end{equation}
    \begin{equation}
        w_i (t+1)= -w_i (t) + 2\cdot r_2\cdot(x_i (t) - g_i (t)) + 2\cdot r_2\cdot(x_i (t) - p_i (t))
    \end{equation}
Where $l_i$ is the best-performing neighbor of the individual particle, thus following multidimensional ring topology \cite{topology} $v_i$ is velocity,whose goal is to pull particles towards the global solution.
    \\
    Secondly - custom coefficients are currently more of the structure element. As of now the main development and movement of the whole colony is more dependent on the other components of our algorithm, the purpose of the individual behaviour of each plant is ensuring of convergence, and allows for individual development of each plant genom. The choice of -2 as coefficients for social and cognitive elements in root factor is based on the fact that the easiest solution would be simply making the movement negative without any randomness. However, we insist on remaining some stochastic properties to ensure the certain level of diversity in colony behaviour. 
    \item \underline{Collective} 
    \\
\item \underline{Theoretical Framework and Application}

In the realm of swarm intelligence algorithms, the utilization of weighted averages as a pivotal mechanism for search region exploitation is not a novel concept. This strategy has been effectively integrated into various algorithms, such as the Gravitational Search Algorithm (GSA)-PSO, Quantum Particle Swarm Optimization (QPSO), among others. Commonly, this approach operates under the hypothesis that, given substantial exploration capabilities or an efficacious initialization strategy, the global optimum tends to be proximal to the swarm's centroid, particularly nearer to the more proficient particles.
\item \underline{Distinctive Aspects of weighted average application}

In our approach, the implementation of weighted averages diverges in two significant facets:

    Historical Data Integration: \\Contrary to the traditional method of considering only the current positional and fitness data of the swarm, our approach incorporates a preserved historical perspective, as delineated in the section 'Preserved History'. This historical integration, while enriching the algorithm with past insights, also introduces potential complexities in heterogeneous search landscapes. Specifically, if the swarm initiates exploration in areas with diverse directional tendencies, discrepancies in historical memories may lead to a decoupling of the swarm's trajectory from its accumulated memory. Notwithstanding these challenges, our methodology exhibits robustness and effectiveness in multimodal functions with stable tendencies, effectively discerning promising search directions.

    Dual-Function Utilization: \\The weighted average method serves a bifurcated purpose in our framework. Primarily, it contributes to generating a constellation of central points (as detailed later), thereby accelerating convergence. Its ancillary role lies in augmenting exploration capabilities by introducing an adaptive component into collective reprojection, which is steered by the backward average. This backward average, illustrated in sections A, B, C, is presumed to exhibit a positive divergence in environmental characteristics, as manifested through the integral of the objective function across their respective domains.//
    Assumption 1
    \begin{equation}
    \int_{F-\epsilon}^{F+\epsilon} g(x) \, dx < \int_{B-\epsilon}^{B+\epsilon} g(x) \, dx
    \end{equation}
    Where epsilon is a distance metric,restricting the environments around the backward and forward weighted averages, B is the backward center, F is the forward weighted center, and g(x) is the function being minimized
    This assumption underpins the generation of subsequent points and collective reprojection in our model.

Adaptive Reprojection Magnitude

A crucial element in our algorithmic design is the incorporation of maximum deviation, which serves to tailor the magnitude of reprojection. This adjustment is attuned to the swarm's behavior and is contingent on the extent of exploitation within the current search region. By modulating reprojection in response to these dynamics, our methodology achieves a nuanced balance between exploration and exploitation, enhancing the overall efficacy of the search process.
    \\ 
    Firstly, here is the general equation of the weighted center. 
    \\
    \begin{equation}
        \mu_w=\frac{\sum\limits_{i=1}^{N_p}( w_i \cdot X_i)}{\sum\limits_{i=1}^{N_p} w_i}
    \end{equation}
    
    where $w_i$ is the weight of $i$-th particle, $X_i$ is a coordinates vector of the particle, and $N_p$ is the number of particle in preserved memory. Here we see the equation of the weight formed by $j$-th function for an $i$-th preserved particle
    
    \begin{equation}
        w_{i,j}=\frac{f_j\left(t_i\right)}{a_i^\alpha}
    \end{equation}
    
    \begin{equation}
        t_i=\frac{c_i-\mu_c}{\sigma_c}
    \end{equation}
    
Let $c_i$ represent the fitness of the $i$-th particle within the swarm. The parameters $\mu_c$ and $\sigma_c$ denote the mean and standard deviation, respectively, of the fitness values as recorded in the preserved memory of the swarm. This historical data provides a contextual backdrop against which current fitness evaluations are assessed.

The term $t_i$ refers to the preliminary weight assigned to the $i$-th particle. This initial weight is subject to further transformation, adhering to a specified mathematical relationship detailed in the equation.\\

An essential aspect of the weighting mechanism is the incorporation of particle age, represented as $a_i^\alpha$ . Here,  $a_i$ denotes the age of the $i$-th particle, and $\alpha$ is a tunable exponent with the constraint of being positive. The parameter $\alpha$  plays a pivotal role in modulating the swarm's sensitivity to temporal aspects of the search process. A higher value of $\alpha$  increases the swarm's responsiveness to solutions discovered in later generations, imbuing the search with a temporal depth. Conversely, reducing $a_i^\alpha$ shifts the focus predominantly towards immediate fitness values, deemphasizing the temporal dimension of the search.

In addition, $f_j(x)$ symbolizes one of the selected functions from a predefined set, each corresponding to a specific mean value within the algorithm's structure. The choice of these functions is fundamentally architectural, restricted to those with monotonically increasing derivatives. In scenarios where minimization is the objective, such as minimizing a model's error, these functions are typically raised to the power of -1. This inversion aligns with the minimization goal, ensuring that the algorithm prioritizes lower function values, thereby steering the search towards optimal solutions.

This intricate weighting scheme, combining fitness evaluation, historical context, particle age, and function selection, exemplifies the algorithm's comprehensive approach. It encapsulates both the immediate and historical performance of particles, thereby balancing exploration and exploitation in the swarm's collective search dynamics.
    \\
In the subsequent stage of the algorithm, particle generation is meticulously orchestrated based on the number and nature of weight functions employed. This process is quantitatively defined as follows:

The algorithm initiates the creation of \( 2 \cdot f_n \) new particles, where \( f_n \) is a user-defined parameter that signifies the count of distinct weight functions being utilized within the algorithm's framework. This parameter, \( f_n \), is central to determining the diversity and quantity of the new particles generated.

For each weight function, the algorithm generates a pair of particles. The first particle in this pair is directly derived from the weighted average calculated using the corresponding weight function. This weighted average particle embodies the collective characteristics influenced by the specific weight function, encapsulating the aggregate information processed through that function.

The second particle in each pair is generated through a process of reprojection. This reprojection is executed from a common reference point, denoted as \( \mu_b \), or backward weighted average. The reprojection involves recalculating the position of the particle by orienting it away from this less favorable average, thereby steering it towards potentially more promising regions of the search space.

This dual particle generation mechanism – one representing the weighted average and the other being reprojected from an unfavorable mean – is a strategic approach. It ensures that for each weight function, there is not only a representation of the collective trends dictated by that function but also an exploration away from less desirable regions. This balance between representation and exploration is crucial in enhancing the robustness and efficiency of the algorithm in navigating complex search spaces.
    \begin{equation}
        \mu_b = \frac{\sum\limits_{i=1}^{N} (t_i\cdot X_i)}{\sum\limits_{i=1}^{N}t_i}
    \end{equation}
    Here is the equation for the bad mean, where N is a size of the swarm, $X_i$ is a coordinate vector of $i$-th particle, and $t_i$ is the inversed normalized fitness value of $i$-th particle(in case of maximization tasks it should the inversion is not required). 
    \\
    Then, the second particle in pair, is generated by reprojection of ,,bad mean’’ across the corresponding weighted average, according to the equation
    
    \begin{equation}
        r = \mu_w- \mu_b
    \end{equation}
    
    \begin{equation} 
        x_r=\mu_b+r\cdot (1+ \frac{\sigma}{\| r\|_n})
    \end{equation}
    
    Where $\mu_b$ is backward weighted average, r is the D-dimensional coordinates vector, containing the difference in position between the corresponding forward weighted average and the backward weighted center, $\|r\|_n$ is equivalent to Euclidean distance between them, $x_r$ is a coordinates vector of a new particle, and  $\sigma$ is a span measure for the search space, considered at the moment.
    
    \begin{equation}
        \sigma = \left[\sigma_1,\sigma_2,\ldots,\sigma_{D-1},\sigma_D\right]_{max}
    \end{equation}
    As we can see here, sigma is an infinity order norm of the D-dimensional vector of standard deviations of  coordinates in preserved history.
    The primary mechanism for directed movement within the swarm in our algorithm is the process of reprojection based on the optimal neighbors of each particle. This method may superficially resemble the numerical gradient estimation technique using finite differences, but fundamental differences exist between the two approaches.

In the finite difference method, a predetermined small value, commonly denoted as epsilon (\(\epsilon\)), is used to approximate the gradient of a function. This method involves calculating the function values at points separated by \(\epsilon\) and using these values to estimate the slope or gradient of the function (REF: [Standard reference on numerical methods]). Essentially, it employs a form of linear approximation where the choice of \(\epsilon\) and the specific points of evaluation determine the accuracy of the gradient estimate.

Contrastingly, our algorithm adopts a distinct strategy where the 'epsilon' in each iteration, conceptualized as the distance between the backward and forward weighted averages, is not just a tool for estimation but serves as a significant metric in itself. This distance reflects the dynamic scope of the search in the problem space and is adaptive to the swarm's behavior.

Moreover, in classical numerical gradient methods, irrespective of the chosen \(\epsilon\), there is always a fixed reference point. The method's nature (forward or backward difference) is defined by how this reference point is utilized. In our approach, however, the concept of a reference point is more fluid. The 'best neighbors' or the particle for which these neighbors are optimal can alternately function as reference points. This flexibility is crucial as it shifts the focus from absolute positions in the search space to the relative positioning and interactions within the swarm. We emphasize the relationships and dynamics between particles rather than their specific coordinates, aligning with the underlying principles of swarm intelligence where collective behavior and interaction patterns are paramount. 

This fundamental difference in approach highlights the novelty of our method. By focusing on relative positions and adaptive distances within the swarm, our algorithm can effectively navigate complex search spaces, leveraging the collective intelligence and adaptability inherent in swarm-based techniques.
\end{enumerate}
\item \underline{Natural selection}
After each iteration, $2\cdot f_n$ particles with the worst fitness value are removed from the swarm – according to the number of new particles generated by collective reprojection, thus leaving the size of the swarm the same.

\end{enumerate}
Here are the steps of ALMI-PSO
\begin{algorithm}
  \caption{ALMI-PSO}
  \begin{algorithmic}
    \For{\text{generation = 1,2,...,t}}:
        \State Update the preserved history
        \State Generate new particles based on preserved history and evaluate them
        \State Natural selection
        \State Update current best solution

        \For{\text{plant = 1,2,...,N}}:
            \State Update personal and social information
            \State Update velocity factor
            \State Update the coordinates
            \State Update fitness values
        \EndFor

    \EndFor
  \end{algorithmic}
\end{algorithm}

\section{Algorithm Comparison on functions}

We are going to assess the performance of our algorithm  on two sets of benchmark functions.Firstly, we evaluate it's performance in comparison to other algorithms on functions suiting it's main abilities, and then we compare the results to the evaluation on CEC2017 test functions.Data about the performance of the other algorithms\cite{XU201933}\cite{1202264}\cite{parsopoulos2019upso}\cite{LIN2022469} is taken from the article \cite{LIN2022469},involving a suitably diverse set of functions for test.

\begin{table}[h]
\caption{Expressions of tendency-possessing functions}
\renewcommand\arraystretch{2.5}
\begin{tabular}[c]{|c|c|c|}
    \hline
    Name of the function & Way to refer & Expression\\ [0.6em]
    \hline Sphere &  STF1 & $\sum\limits_{i=1}^{M}{x_i}^2$\\ 
    \hline Schwefel’s P2.22 & STF2 & $\sum_{i=1}^{n} \sum_{j=1}^{i} x_j$\\ 
    \hline Quadric & STF3 & $\sum_{i=1}^{n} |x_i| + \prod_{i=1}^{n} |x_i|$\\  
    \hline Step & STF4 & $\sum_{i=1}^{n} \left\lfloor x_i + 0.5 \right\rfloor^2 $\\ 
    \hline Schwefel's & STF5& $\sum_{i=1}^{n} random[0,1]  \cdot x_i^4$ \\ 
    \hline Rastrigrin & STF6 &$10n + \sum_{i=1}^{n} ( x_i^2 - 10 \cos(2\pi x_i)$ \\
    \hline Ackley & STF7& $-20 \exp\left(-0.2 \sqrt{\frac{\sum_{i=1}^{n} x_i^2}{n}}\right) - \exp\left(\frac{\sum_{i=1}^{n} \cos(2\pi x_i)}{n}\right) + 20 + \exp(1)$ \\ 
    \hline Griewank & STF8 &$f_{10}(x) = 1 + \frac{\sum_{i=1}^{n} x_i^2}{4000} - \prod_{i=1}^{n} \cos\left(\frac{x_i}{\sqrt{i}}\right)$\\ 
    \hline
\end{tabular}
\end{table}

\begin{table}[h!]
\caption{Bounds and optimal solutions of manipulating variables}
\renewcommand\arraystretch{1.5}
\begin{tabular}[c]{|c|c|c|}
    \hline
    Function & Initialization variable bounds & Optimum\\
    \hline $STF1$ & $[-100, 100]$M & 0\\
    \hline $STF2$ & $[-10, 10]$M & 0\\
    \hline $STF3$ & $[-100, 100]$M & 0 \\
    \hline $STF4$ & $[-100, 100]$M & 0 \\
    \hline $STF5$ & $[-500, 500]$M & 0 \\
    \hline $STF6$ & $[-5.12, 5.12]$M& 0 \\
    \hline $STF7$ & $[-32, 32]$M & 0 \\
    \hline $STF8$ & $[-600, 600]$M & 0 \\
    \hline
\end{tabular}
\end{table}

In our case, we use two function groups-polynomial and exponential ${{(e}^x)}^{-1}, {{(x}^3)}^{-1}$ – as a basic example of the exponential function and the cubic polynomial, to emphasize the importance of the power being odd in case of polynomial space. These are the general functions, which were chosen to show the general performance, without tuning the architecture to suit the certain purpose better. Raising them to the power of -1 directly corresponds to the fact that we will be solving minimization problems.
\\
\indent The age is going to have the power of 1 to preserve generality, colony size to 10, so that one-fifth of particles would be the generated ones, memory length to 60 and maximum number of function evaluations to 10000*D, where D is the dimensionality, making 30 runs for more thorough analysis.
\\

\begin{enumerate}
    \item \textbf{Analysis with Functions Exhibiting Stable Tendencies:} The initial phase of our testing regime involves a set of functions characterized by stable tendencies. These functions are selected to assess the algorithm's ability to efficiently navigate optimization landscapes with consistent directional gradients. This evaluation will provide insights into the algorithm's core optimization proficiency in relatively straightforward scenarios.

    \item \textbf{Assessment Using Rotated Function Sets:} Subsequently, the algorithm will be subjected to a series of rotated functions. This phase is critical for demonstrating the algorithm's invariance to rotational transformations of the search space. By evaluating the algorithm's performance on these rotated functions, we aim to ascertain its independence from orientation biases, a vital attribute for effective multidimensional optimization.

    \item \textbf{Testing on Functions with Diverse Directional Tendencies:} The final set of benchmark tests comprises functions that deviate from the aforementioned stable tendencies. These functions are characterized by regions exhibiting diverse directional tendencies, potentially leading to discrepancies in the algorithm's historical memory. Such a testing environment is crucial for examining the resilience of the algorithm when confronted with complex and unpredictable search landscapes. The presence of varied directional cues challenges the algorithm's ability to maintain a coherent search trajectory, particularly in relation to its accumulated historical data.
\end{enumerate}

For each category of benchmark tests, we plan to conduct a thorough comparative analysis of the algorithm's performance. This will involve not only direct comparisons within each function set but also cross-comparisons between the different sets to gain a comprehensive understanding of the algorithm's versatility and adaptability across diverse optimization challenges.

Furthermore, to quantitatively assess the diversity and differentiation among the performance of various algorithms, including ours, we will employ the Friedman test. This non-parametric statistical test will be applied to each group of functions, provided they exhibit sufficient variability in performance. The Friedman test is particularly suitable for identifying instances where different algorithms may exhibit dominance under specific conditions or within certain function categories. This statistical analysis will add a layer of rigor to our evaluation, ensuring that our findings are grounded in statistically valid comparisons.

In summary, this multi-faceted testing approach is designed to provide a holistic assessment of the algorithm's capabilities, revealing strengths, potential limitations, and areas for further refinement. The outcome of these tests will offer valuable insights, guiding future developments and optimizations of the algorithm.
Here, we are going to define the variable O- any number closer to the value of a function in a global optima than a certain ,,accepted distance’’ will be simply set to O,and considered global optimum of a function. As accepted distance, we use $4.94\cdot {10}^{-324}$, as we operate in 64bit numbers, and any distance lower is not deemed perceivable.

\begin{table}[h!]
    \caption{Test results on stable-tendency functions - STF1}
    \renewcommand\arraystretch{0.6}
    \centering
    \hspace*{-2.5cm}\begin{tabular}{|cc|cc|cc|cc|cc|cc|cc|cc|}
        \hline 
        Statistic&ALMI-PSO & FDR-PSO & UPSO & TSL-PSO & GGL-PSOD & DMSDL-PSO & MLDE-PSO \\
        \hline
        FV &0.000E + 00 &2.182E-30 & 1.283E-149 & 1.283E-149 & 2.557E-203 & 1.959E-57 & 0.000E + 00 \\\hline
        SP &3.63E + 01& 1.231E + 05 & 1.148E + 05 & 2.333E + 04 & 3.313E + 04 & 1.493E + 05 & 9.012E + 04 \\\hline
        SR & 100.00\%& 100.00\% & 100.00\% & 100.00\% & 100.00\% & 100.00\% & 100.00\% \\\hline
        Rank &1 & 3 & 4 & 2 & 6 & 5 & 2 \\
        \hline
    \end{tabular}
\end{table}

\begin{table}[h!]
    \caption{Test results on stable-tendency functions - STF2}
    \renewcommand\arraystretch{1.0}
    \centering
    \hspace*{-2.5cm}\begin{tabular}{|c|c|c|c|c|c|c|c|}
        \hline 
        Statistic& ALMI-PSO & FDR-PSO & UPSO & TSL-PSO & GGL-PSOD & DMSDL-PSO & MLDE-PSO \\
        \hline
        BF 2  &0.000E + 00 & 3.696E-63 & 3.696E-63 & 1.490E-97 & 5.398E-33 & 1.208E-34 & 0.000E + 00 \\\hline
        SP &0.000E + 01 &1.225E + 05 &3.136E + 04 &4.434E + 04 &1.657E + 05 &7.310E + 04 &2.663E + 04 \\\hline
        SR & 100.00\% & 100.00\% & 100.00\% & 100.00\% & 100.00\% & 100.00\% & 100.00\% \\\hline
        Rank &1 & 6 & 3 & 4 & 2 & 5 & 2 \\
        \hline
    \end{tabular}
\end{table}

\begin{table}[h!]
    \caption{Test results on stable-tendency functions - STF3}
    \renewcommand\arraystretch{1.0}
    \centering
    \hspace*{-2.5cm}\begin{tabular}{|c|c|c|c|c|c|c|c|}
        \hline 
        Statistic&ALMI-PSO & FDR-PSO & UPSO & TSL-PSO & GGL-PSOD & DMSDL-PSO & MLDE-PSO \\
        \hline
        BF 3 & 0.000E + 00 & 4.603E-150 & 4.603E-150 & 2.000E-203 & 8.137E-56 & 5.657E-65 & 0.000E + 00 \\\hline
        SP &4.4E +01 & 1.185E + 05 & 2.588E + 04 & 3.329E + 04 & 1.543E + 05 & 5.510E + 04 & 2.009E + 04 \\\hline
        SR  &100.00\%& 100.00\% & 100.00\% & 100.00\% & 100.00\% & 100.00\% & 100.00\% \\\hline
        Rank &1 & 3 & 4 & 2 & 6 & 5 & 2 \\
        \hline
    \end{tabular}
\end{table}

\begin{table}[h!]
    \caption{Test results on stable-tendency functions - STF4}
    \renewcommand\arraystretch{1.0}
    \centering
    \hspace*{-2.5cm}\begin{tabular}{|c|c|c|c|c|c|c|c|}
        \hline 
        Statistic&ALMI-PSO & FDR-PSO & UPSO & TSL-PSO & GGL-PSOD & DMSDL-PSO & MLDE-PSO \\
        \hline
        FV & 0.000E + 00& 1.667E-01 & 1.667E-01 & 0.000E + 00 & 0.000E + 00 & 0.000E + 00 & 0.000E + 00 \\\hline
        SP  & 2.8E + 01& 1.474E + 05 & 1.138E + 04 & 1.652E + 04 & 1.115E + 05 & 2.296E + 04  & 7.964E + 03 \\\hline
        SR & 100.00\% & 83.33\% & 100.00\% & 100.00\% & 100.00\% & 100.00\% & 100.00\% \\\hline
        Rank & 1 & 3 & 3 & 2 & 2 & 2 \\
        \hline
    \end{tabular}
\end{table}

\begin{table}[h!]
    \caption{Test results on stable-tendency functions - STF5}
    \renewcommand\arraystretch{1.0}
    \centering
    \hspace*{-2.5cm}\begin{tabular}{|c|c|c|c|c|c|c|c|}
        \hline 
        Statistic &ALMI-PSO& FDR-PSO & UPSO & TSL-PSO & GGL-PSOD & DMSDL-PSO & MLDE-PSO \\
        \hline
        FV & 0.000E + 00& 2.532E-03 & 2.532E-03 & 3.517E-03 & 2.471E-03 & 2.094E-03 & 1.953E-04 \\\hline
        SP & 4.0E + 01& Inf & Inf & Inf & Inf & Inf & Inf \\\hline
        SR & 100.00\%& 0.00\% & 0.00\% & 0.00\% & 0.00\% & 0.00\% & 0.00\% \\\hline
        Rank &1 & 4 & 5 & 6 & 3 & 2 & 7 \\
        \hline
    \end{tabular}
\end{table}

\begin{table}[h!]
    \caption{Test results on stable-tendency functions - STF6}
    \renewcommand\arraystretch{1.0}
    \centering
    \hspace*{-2.5cm}\begin{tabular}{|c|c|c|c|c|c|c|c|}
        \hline 
        Statistic &ALMI-PSO& FDR-PSO & UPSO & TSL-PSO & GGL-PSOD & DMSDL-PSO & MLDE-PSO\\
        \hline
        FV & 0.000E + 00&2.802E + 01 & 2.802E + 01 & 0.000E + 00 & 5.329E-16 &  0.000E + 00 & 0.000E + 00 \\\hline
        SP & 4.0E + 01& Inf & Inf & 3.349E + 04 & 2.178E + 05 & 6.946E + 04 & 2.322E + 04 \\\hline
        SR & 100.00\%& 0.00\% & 0.00\% & 100.00\% & 100.00\% & 100.00\%  & 100.00\% \\\hline
        Rank &1 & 3 & 4 & 1 & 2 & 1  & 1 \\
        \hline
    \end{tabular}
\end{table}

\begin{table}[h!]
    \caption{Test results on stable-tendency functions - STF7}
    \renewcommand\arraystretch{1.0}
    \centering
    \hspace*{-2.5cm}\begin{tabular}{|c|c|c|c|c|c|c|c|}
        \hline 
        Statistic &ALMI-PSO & FDR-PSO & UPSO & TSL-PSO & GGL-PSOD & DMSDL-PSO & MLDE-PSO \\
        \hline
        FV & 0.000E + 00 & 2.315E-14 & 2.315E-14 & 3.417E-14 & 8.112E-15 & 6.454E-15 &  8.882E-16 \\\hline
        SP & 2.002E + 02 & 1.300E + 05 & 4.545E + 04 & 4.354E + 04 & 1.688E + 05 & 7.549E + 04 & 2.658E + 04 \\\hline
        SR & 100.00\%& 96.67\% & 100.00\% & 100.00\% & 100.00\% & 100.00\% & 100.00\% \\\hline
        Rank &1 & 5 & 6 & 7 & 4 & 3 & 2 \\
        \hline
    \end{tabular}
\end{table}

\begin{table}[h!]
    \caption{Test results on stable-tendency functions - STF8}
    \renewcommand\arraystretch{1.0}
    \centering
    \hspace*{-2.5cm}\begin{tabular}{|c|c|c|c|c|c|c|c|}
        \hline 
        Statistic&ALMI-PSO & FDR-PSO & UPSO & TSL-PSO & GGL-PSOD & DMSDL-PSO & MLDE-PSO \\
        \hline
        FV & 0.000E + 00 & 1.532E-02 & 1.532E-02 & 6.049E-13 & 3.281E-03 &  0.000E + 00 & 0.000E + 00 \\\hline
        SP & 4.0E + 01 & 6.346E + 05 & 1.160E + 05 & 3.443E + 04 & 2.423E + 05 & 6.193E + 04  & 1.901E + 04 \\\hline
        SR & 100.00\% & 36.67\% & 80.00\% & 100.00\% & 76.67\% & 100.00\% & 100.00\% \\\hline
        Rank &1 & 5 & 6 & 3 & 4 & 2 & 2  \\
        \hline
    \end{tabular}
\end{table}
From these results we can see the superiority of the ALMI-PSO on tendency-possessing functions both in term of success rate and convergency speed

\begin{table}[h!]
    \caption{Test results on rotated functions - RF1}
    \renewcommand\arraystretch{1.0}
    \centering
    \hspace*{-2.5cm}\begin{tabular}{|c|c|c|c|c|c|c|c|}
        \hline 
        Statistic & ALMI-PSO & FDR-PSO & UPSO & TSL-PSO & GGL-PSOD & DMSDL-PSO & MLDE-PSO \\
        \hline
        FV & 0.000E + 00 & 2.156E-114 & 1.252E-98 & 2.340E-121 & 1.161E-52 & 1.631E-43 & 5.373E-28 \\\hline
        SP & 3.63E + 01 & 1.198E + 05 & 3.040E + 04 & 4.859E + 04 & 1.535E + 05 & 7.736E + 03 & 1.009E + 05 \\\hline
        SR & 100.00\% & 100.00\% & 100.00\% & 100.00\% & 100.00\% & 100.00\% & 100.00\% \\\hline
        Rank & 1 & 3 & 4 & 2 & 5 & 6 & 7 \\
        \hline
    \end{tabular}
\end{table}

\begin{table}[h!]
    \caption{Test results on rotated functions - RF2}
    \renewcommand\arraystretch{1.0}
    \centering
    \hspace*{-2.5cm}\begin{tabular}{|c|c|c|c|c|c|c|c|}
        \hline 
        Statistic & ALMI-PSO & FDR-PSO & UPSO & TSL-PSO & GGL-PSOD & DMSDL-PSO & MLDE-PSO \\
        \hline
        BF 2 & 0.000E + 00 & 3.911E-01 & 1.094E + 01 & 1.123E + 01 & 4.965E-04 & 1.991E-05 & 3.290E-11 \\\hline
        SP & 4.2E + 01 & Inf & Inf & Inf & 2.003E + 05 & 1.458E + 06 & 2.017E + 05 \\\hline
        SR & 100.00\% & 0.00\% & 0.00\% & 0.00\% & 93.33\% & 20.00\% & 100.00\% \\\hline
        Rank & 1 & 5 & 6 & 7 & 4 & 3 & 2 \\
        \hline
    \end{tabular}
\end{table}

\begin{table}[h!]
    \caption{Test results on rotated functions - RF3}
    \renewcommand\arraystretch{1.0}
    \centering
    \hspace*{-2.5cm}\begin{tabular}{|c|c|c|c|c|c|c|c|}
        \hline 
        Statistic & ALMI-PSO & FDR-PSO & UPSO & TSL-PSO & GGL-PSOD & DMSDL-PSO & MLDE-PSO \\
        \hline
        FV & 0.000E + 00 & 6.727E-19 & 1.288E-10 & 1.872E-31 & 3.761E-12 & 8.727E-12 & 5.154E-16 \\\hline
        SP & 4.4E + 01 & 1.919E + 05 & 1.944E + 05 & 8.210E + 04 & 2.307E + 05 & 2.497E + 05 & 1.568E + 05 \\\hline
        SR & 100.00\% & 100.00\% & 100.00\% & 100.00\% & 100.00\% & 100.00\% & 100.00\% \\\hline
        Rank & 1 & 3 & 7 & 2 & 5 & 6 & 4 \\
        \hline
    \end{tabular}
\end{table}

\begin{table}[h!]
    \caption{Test results on rotated functions - RF4}
    \renewcommand\arraystretch{1.0}
    \centering
    \hspace*{-2.5cm}\begin{tabular}{|c|c|c|c|c|c|c|c|}
        \hline 
        Statistic & ALMI-PSO & FDR-PSO & UPSO & TSL-PSO & GGL-PSOD & DMSDL-PSO & MLDE-PSO \\
        \hline
        FV & 0.000E + 00 & 4.100E + 00 & 0.000E + 00 & 9.667E-01 & 3.333E-02 & 0.000E + 00 & 0.000E + 00 \\\hline
        SP & 4.4E+01 & Inf & 3.567E + 04 & 2.204E + 05 & 1.269E + 05 & 3.783E + 04 & 4.726E + 04 \\\hline
        SR & 100.00\% & 0.00\% & 100.00\% & 96.67\% & 96.67\% & 100.00\% & 100.00\% \\\hline
        Rank & 1 & 5 & 2 & 4 & 3 & 2 & 2 \\
        \hline
    \end{tabular}
\end{table}

\begin{table}[h!]
    \caption{Test results on rotated functions - RF5}
    \renewcommand\arraystretch{1.0}
    \centering
    \hspace*{-2.5cm}\begin{tabular}{|c|c|c|c|c|c|c|c|}
        \hline 
        Statistic & ALMI-PSO & FDR-PSO & UPSO & TSL-PSO & GGL-PSOD & DMSDL-PSO & MLDE-PSO \\
        \hline
        FV & 0.000E + 00 & 2.637E-03 & 1.347E-02 & 5.054E-03 & 2.128E-03 & 3.495E-03 & 5.338E-03 \\\hline
        SP & 4.0E+01 & Inf & Inf & Inf & Inf & Inf & Inf \\\hline
        SR & 100.00\% & 0.00\% & 0.00\% & 0.00\% & 0.00\% & 0.00\% & 90.00\% \\\hline
        Rank & 1 & 3 & 7 & 5 & 2 & 4 & 6 \\
        \hline
    \end{tabular}
\end{table}

\begin{table}[h!]
    \caption{Test results on rotated functions - RF6}
    \renewcommand\arraystretch{1.0}
    \centering
    \hspace*{-2.5cm}\begin{tabular}{|c|c|c|c|c|c|c|c|}
        \hline 
        Statistic & ALMI-PSO & FDR-PSO & UPSO & TSL-PSO & GGL-PSOD & DMSDL-PSO & MLDE-PSO \\
        \hline
        FV & 0.000E + 00 & 4.985E + 01 & 9.170E + 01 & 8.626E + 01 & 3.323E + 01 & 3.293E + 01 & 1.810E + 02 \\\hline
        SP & 4.0E+01 & Inf & Inf & Inf & Inf & Inf & 3.442E + 04 \\\hline
        SR & 100.00\% & 0.00\% & 0.00\% & 0.00\% & 0.00\% & 0.00\% & 100.00\% \\\hline
        Rank & 1 & 4 & 6 & 5 & 3 & 2 & 8 \\
        \hline
    \end{tabular}
\end{table}

\begin{table}[h!]
    \caption{Test results on rotated functions - RF7}
    \renewcommand\arraystretch{1.0}
    \centering
    \hspace*{-2.5cm}\begin{tabular}{|c|c|c|c|c|c|c|c|}
        \hline 
        Statistic & ALMI-PSO & FDR-PSO & UPSO & TSL-PSO & GGL-PSOD & DMSDL-PSO & MLDE-PSO \\
        \hline
        FV & 0.000E + 00 & 1.866E + 00 & 2.281E + 00 & 2.519E + 00 & 1.036E-14 & 7.832E-09 & 8.882E-16 \\\hline
        SP & 0.000E + 00 & 4.340E + 06 & Inf & 7.666E + 05 & 1.746E + 05 & 1.833E + 05  & 2.795E + 04 \\\hline
        SR & 0.00\% & 6.67\% & 0.00\% & 30.00\% & 100.00\% & 100.00\% & 100.00\% \\\hline
        Rank & 1 & 5 & 7 & 8 & 2 & 4  & 1 \\
        \hline
    \end{tabular}
\end{table}

\begin{table}[h!]
    \caption{Test results on rotated functions - RF8}
    \renewcommand\arraystretch{1.0}
    \centering
    \hspace*{-2.5cm}\begin{tabular}{|c|c|c|c|c|c|c|c|}
        \hline 
        Statistic & ALMI-PSO & FDR-PSO & UPSO & TSL-PSO & GGL-PSOD & DMSDL-PSO & MLDE-PSO \\
        \hline
        FV & 0.000E + 00 & 8.865E-03 & 2.712E-03 & 6.271E-04 & 3.122E-03 & 4.138E-12 & 0.000E + 00 \\\hline
        SP & 0.000E + 00 & 6.526E + 05 & 1.835E + 05 & 5.821E + 04 & 3.157E + 05 & 1.358E + 05  & 2.145E + 04 \\\hline
        SR & 100\% & 36.67\% & 70.00\% & 100.00\% & 66.67\% & 100.00\% & 100.00\% \\\hline
        Rank & 1 & 8 & 6 & 5 & 7 & 3  & 1 \\
        \hline
    \end{tabular}
\end{table}
Here we can accept the rotation-invariance abilities of ALMI-PSO.\\

\begin{table}[h!]
    \caption{Test results on CEC2017 functions}
    \renewcommand\arraystretch{1.0}
    \centering
    \hspace*{-2.5cm}\begin{tabular}{|c|c|c|c|c|c|c|c|}
        \hline 
        Function&ALMI-PSO & FDR-PSO & UPSO & TSL-PSO & GGL-PSOD & DMSDL-PSO & MLDE-PSO \\
        \hline
        F1 
 & 5.71e+08 & 1.677E + 03 &3.377E + 03&1.677E + 03& 2.876E + 03 & 7.447E + 01  & 0.000E + 00 \\\hline
        F2  & 3.92e+03& 4.914E + 10 & 1.486E + 16 & 4.914E + 10 & 1.509E + 06 & 1.948E + 10  & 1.667E + 00 \\\hline
        F3 & 3.72e+02 & 4.145E-08 & 4.754E + 01 & 4.145E-08 & 1.103E-07 & 1.199E-05  & 2.519E-09 \\\hline
        F4 &8.62e+02 & 2.465E + 01 & 1.240E + 02 & 2.465E + 01 & 6.404E + 01 & 1.728E + 00 & 4.913E + 01 \\\hline
        F5 &5.00e+02 & 5.877E + 01 & 8.236E + 01 & 5.877E + 01 & 3.042E + 01 & 3.841E + 01  & 2.584E + 01 \\\hline
        F6 &6.31e+03 & 1.860E-01 & 1.294E + 00 & 1.860E-01 & 5.917E-08 &  1.137E-13 & 1.091E-07 \\\hline
        F7 &7.49e+02 & 1.009E + 02 & 1.334E + 02 & 1.009E + 02 & 6.282E + 01 & 7.200E + 01  & 6.198E + 01 \\\hline
        F8 &8.01e+02 & 5.860E + 01 & 8.430E + 01 & 5.860E + 01 & 3.279E + 01 & 4.059E + 01  & 2.790E + 01 \\
        \hline
    \end{tabular}
\end{table}
And here,the group with negative results is presented,where there is a drastic performance difference seen.
\\ \\ \\ \\ \\ \\ \\ \\ \section{Statistical Analysis}
Analysis of Diversity Within Groups

In our research, we employed the Friedman test to statistically evaluate the performance variations of seven distinct implementations of Particle Swarm Optimization (PSO) across eight functions. This test examines the null hypothesis that repeated measurements on the same entities follow an identical distribution. It's noteworthy that the test's reliability, predicated on a chi-squared distribution assumption, is most robust when the sample size (n) exceeds 10 and includes more than six repeated measurements.

For the number of function evaluations (SP), a significant diversity in outcomes was evident at the 0.01 significance level, contingent on the algorithm variant utilized, keeping our requirement of diverse algorithm behaviours intact, in order to keep different efficiency levels and provide more unbiased analysis. Analogously, the function values (FV) demonstrated a similar pattern of significant variation at the same statistical threshold. These findings are particularly salient in the context of the ALMI-PSO algorithm, which consistently achieved the top rank (rank 1) across evaluations. This prominence was especially pronounced for rotated functions and stable tendency functions. However, in the case of the cec2017 functions, although a significant difference at the 1
Cross-Group Comparative Analysis
Wilcoxon Signed-Rank Test for Comparative Assessment

Further analysis was conducted using the Wilcoxon signed-rank test to compare paired samples' distributions. This test was applied to evaluate the consistency of rankings between different types of functions. It was observed that for rotated functions and stable tendency functions, the ALMI-PSO algorithm maintained consistent rankings. Conversely, for cec2017 functions, the rankings differed significantly at the 1
Statistical Results

The results for RBF and stable tendency functions (both for FV and SP) indicated significant differences, as evidenced by the obtained p-values:

\begin{table}[h!]
\centering
\caption{Friedman Test Results for Function Values (FV)}
\label{tab:stat_results}
\begin{tabular}{|l|l|l|}
\hline
\textbf{Function set} & \textbf{Statistic} & \textbf{p-value} \\ \hline
RF – FV              & 18.857             & 0.0044           \\ \hline
STF – FV              & 29.017             & 6.04e-05         \\ \hline
CEC2017 – FV          & 31.636             & 1.92e-05         \\ \hline
\end{tabular}
\end{table}

\begin{table}[h!]
\centering
\caption{Statistical Test Results for Comparisons}
\label{tab:comparison_results}
\begin{tabular}{|l|l|l|}
\hline
\textbf{Comparison}     & \textbf{Statistic} & \textbf{p-value} \\ \hline
RF vs. CEC2017         & 0.0                & 0.0078           \\ \hline
STF vs. CEC2017         & 0.0                & 0.0078           \\ \hline
\end{tabular}
\end{table}
There was no difference in ranks between RF and STF set,as they were identical

Summary of Findings

The statistical analysis indicates that the ALMI-PSO algorithm demonstrates superior performance for rotated and stable tendency functions. However, its efficacy varies notably with cec2017 functions. These results suggest that ALMI-PSO is particularly effective for specific types of functions, underscoring the importance of algorithm selection based on the function architecture in optimization tasks.

\pagebreak
\hfill \break
\hfill \break
\\
\\
\\
\\
\indent As we can see, in case of functions fulfilling it's assumptions ALMI-PSO consistently outperforms the other versions by a wide margin,not only being more efficient in finding the global minima,but also doing it considerably faster than the other versions, although it's performance drastically falls back on functions falling out of it's range.
\section{Statistical analysis of the results}
To compare the performance results of our algorithm between the functions that contain the tendency trait, and these that do not,we apply the Friedman's test,that shows the test statistic of 2.714

\section{Results on exponential smoothing}
In this section we apply ALMI-PSO to the parameter search of exponential smoothing algorithm 
Instead of doing a costly grid search with computational compexity of \begin{equation}
    {O}(a^{-d})
\end{equation}, meaning exponential scaling in terms of dimensionality and polynomially(of order at least d) in terms of the precision, our algorithm has only quadratic computational complexity in terms of dimensionality, and is constant in terms of precision(if the IEEE standard of floating point precision stays the same, otherwise linearly). 
As evaluation metric we will use MAPE and RMSE on forecasting with different time windows.
Here we display the loss function on different window lengths for a sinusoid sequence:

\subsection{Single Exponential Smoothing}
As we have mentioned before, single exponential smoothing is the baseline version of the algorithm for data without a trend. As a sequence of interest we will incorporate a sinusoid with added gaussian noise
Here are the results:
Differentiable

\subsection{Double Exponential Smoothing}
Double exponential smoothing is an extension of the single exponential smoothing which allows to count trend into a model, while also doubling the number of parameters to search through
Here are the results on different window lengths
\subsection{Triple Exponential Smoothing}
Triple exponential smoothing is another extension, which allows to count in the cyclical patterns of the data, while also adding another parameter to determine
\section{Conclusions}

In this paper, a new version of PSO for determining the parameters of low pass filters, such as exponential smoothing is proposed. By incorporating the learning of the search space via preserving positive experiences, it also speeds up the convergence and offers high capabilities in finding global optima on multimodal functions with stable tendencies towards global optima,although losing in performance on function without that trait. Tests on it's primary application, searching for the perfect smoothing parameter show it's effectiveness.

\section{Declarations}
Ethical Approval: Not applicable
\\Funding: Not applicable 
\\Availability of data and materials : Not applicable
\bibliographystyle{plain}
\bibliography{ALMI_PSO}
\end{document}